\documentclass[11pt]{amsart}

\usepackage{amssymb,latexsym} \usepackage{palatino}
\usepackage{mathrsfs} \input{sommers.tex}

\newcommand{\nilrad}{\mathfrak{n}} 
 \newcommand{\g}{\mathfrak{g}}
\newcommand{\bo}{\mathfrak{b}} \newcommand{\h}{\mathfrak{h}}

\newcommand{\sz}[1]{\|#1\|}
\newcommand{\complex}{\mathbf C}
\newcommand{\supp}{\text{supp}}
\newcommand{\chamber}{\mathfrak{C}}

\newcommand{\posroots}{\ro^+}

\newcommand{\orbit}{\mathcal O}
\newcommand{\sign}{\mathcal S}
\newcommand{\al}{{\alpha}}

\newcommand{\coal}{{\alpha^{\scriptscriptstyle\vee}}}
\newcommand{\ro}{{\Phi}}
\newcommand{\co}{{\Phi^{\scriptscriptstyle\vee}}}

\title{A characterization of Dynkin elements}

\author{Paul E. Gunnells} \address{University of
Massachusetts---Amherst\\ Amherst, MA 01003}
\email{gunnells@math.umass.edu}
\author{Eric Sommers} 
\address{University of
Massachusetts---Amherst\\ Amherst, MA 01003;
Institute for Advanced Study \\
Princeton, NJ 08540}
\email{esommers@math.umass.edu}

\date{20 January 2003}

\begin{document}

\begin{abstract}

We give a characterization of the Dynkin elements of a simple Lie algebra.
Namely, we prove that 
one-half of a Dynkin element is the unique point of minimal length
in its $N$-region.
In type $A_n$ this translates into a statement about the regions
determined by the canonical left Kazhdan-Lusztig cells, which leads to some
conjectures in representation theory.
\end{abstract}

\maketitle

\section{Introduction}

Let $G$ be a connected simple algebraic group over the complex numbers and 
$\g$ its Lie algebra.
We begin by recalling 
the Dynkin-Kostant classification of nilpotent $G$-orbits in 
$\g$ (see \cite{col-mcg:nilp}).

Let $H= [\begin{smallmatrix} 1 & 0  \\
    0 & -1
\end{smallmatrix}]$ and 
$E  = [\begin{smallmatrix} 0 & 1  \\
    0 & 0
\end{smallmatrix}]$  be elements of ${\mathfrak sl}_2(\complex)$. 
%Let $\phi: {\mathfrak sl}_2(\complex) \to \g$ be a homomorphism of Lie algebras.
Let ${\mathcal A}_{hom}$ denote the $G$-conjugacy classes of Lie algebra homomorphisms 
from ${\mathfrak sl}_2(\complex)$ to $\g$. 
Then the map $\Omega: {\mathcal A}_{hom} \to \{ \text{nilpotent $G$-orbits in } \g \}$
given by $\Omega(\phi) = \Ad(G) \phi(E)$
is a bijection:  surjectivity is just the Jacobson-Morozov theorem
and injectivity follows from a theorem of Kostant.
Let $\Upsilon: {\mathcal A}_{hom} \to \{ \text{semi-simple $G$-orbits in } \g \}$
be the map $\Upsilon(\phi) = \Ad(G) \phi(H)$.
A theorem of Mal'cev shows that $\Upsilon$ is injective.
Hence nilpotent orbits in $\g$ are parametrized by the image of $\Upsilon$.  This
finite set was completely determined by Dynkin \cite{dynkin}.

If $\orbit$ is a nilpotent orbit, then an element in the conjugacy class 
$\Upsilon(\Omega^{-1}(\orbit))$ is called a \emph{ Dynkin element} for 
$\orbit$ (or for any element $e \in \orbit$).
If we fix a Cartan subalgebra $\h$ of $\g$ and a set
of simple roots $\Pi \subset \h^*$, then there is a unique Dynkin element $h$
for $\orbit$ which lies in $\h$ and such that $\al(h) \geq 0$ for all $\al \in \Pi$.
We often refer to this element as {\bf the} Dynkin element for $\orbit$;
the numbers $\al(h)$ are assigned to the Dynkin diagram of $G$ and yield
the weighted Dynkin diagram of $\orbit$.  To summarize, this diagram completely
determines the Dynkin elements for $\orbit$ and hence $\orbit$ itself.

The aim of this paper is to give a new characterization of the 
Dynkin elements.  
We begin by partitioning the dominant Weyl chamber into regions indexed by
nilpotent orbits (which we call $N$-regions).   These regions arise by generalizing in the
most straight-forward fashion the notion of $S$-cells 
(also called geometric cells) from the finite Weyl group 
to the affine Weyl group;
they are closely related to Lusztig's $\tilde{S}$-cells \cite{lusztig:affine_spring}.  
Our main theorem states that the point of minimal length in
the closure of the $N$-region corresponding to an orbit $\orbit$ is 
one-half the Dynkin element for $\orbit$.  
In type $A_n$, the closure of an $N$-region coincides 
with the region determined by a canonical left 
Kazhdan-Lusztig cell, so our theorem becomes a statement about Kazhdan-Lusztig cells.

%fix this up a bit
The Dynkin elements for $G$ are important for the representation theory of
an algebraic group of type dual to $G$ over a local field.  
This is a part of
conjectures by Arthur as explained in Vogan \cite{vogan:dynkin}.
Arthur's conjectures posit that one-half of a 
Dynkin element for a group of type dual to the type of $G$
gives rise to a spherical representation of $G$ which is unitary.
This has been proved in some cases:  see  \cite{vogan:dynkin} for references.
In particular, when the local field is the complex numbers, the central role played by one-half
of a Dynkin element in the unitarity question for complex Lie groups 
is one of the key ideas 
in the work of Barbasch and Vogan \cite{bv:unip}. It leads them to introduce the notion
of special unipotent representations for these groups.  
On the algebraic side,
McGovern has shown in the classical groups that these elements yield completely prime primitive
ideals in the universal enveloping algebra of the Lie algebra of type dual to $\g$
(see the last section for a precise statement) \cite{mcgovern:memoirs}.
The main theorem of this paper suggests a generalization about where to look 
for other semi-simple conjugacy classes in $\g$ which may have this property.
We suspect these elements will be important in representation theory and may in particular be 
useful for defining unipotent representations for complex Lie groups.

Here is a brief outline of the paper.  In section 2, we introduce definitions
and state the main theorem, Theorem \ref{main_theorem}.  
In section 3, we give a proof.  In section 4, we make a definition of a property
relating Dynkin elements to $B$-stable ideals in the nilradical of a Borel subalgebra.
We use this definition to explore further properties of the Dynkin elements as minimal
points of certain convex regions.
In the final section, we speculate on how generalizing the idea of the main theorem
to Kazhdan-Lusztig cells and $\tilde{S}$-cells for $G$ of general type may lead to 
an important new class of infinitesimal characters.

We thank J. Humphreys, G. Lusztig, W. McGovern, I. Mirkovic, A. Ram, D. A. 
Vogan, Jr., and the referee
for helpful discussions and/or suggestions.  
We gratefully acknowledge the support of NSF grants DMS-0245580 (P.G.) 
and NSF grants DMS-0201826 and DMS-9729992 (E.S.).

\section{Main Theorem}

Let $(\ro, X, \co, Y)$ be the root datum of a connected
simple algebraic group $G$ defined over the complex numbers
with respect to a maximal torus $T$.  
Let $W$ be the Weyl group.
Fix a Borel subgroup $B$ containing $T$ and let $\ro^+$ 
be the determined positive roots and $\Pi$ the simple roots.  
Let $\g, \bo, \h$ be the Lie algebras of $G$, $B$, $T$, respectively. 
Let $V = Y \otimes \mathbf{R}$.
We make $V$ into a Euclidean space by fixing a $W$-invariant 
positive-definite inner product $( \ , \ )$, with the induced
norm denoted $\sz{\cdot}$.
The inner product is unique up to a non-zero scalar, 
since $V$ is an irreducible representation of $W$.
The closure of the dominant Weyl chamber in $V$ is 
$$\chamber:= \{ x \in V \ | \ \langle \al , x \rangle \geq 0 \text{ for all } \al \in \ro^+\},$$
where $\langle \ , \ \rangle$ is the extension of the 
pairing of weights and coweights.

Inspired by the work of Shi \cite{shi:book}, \cite{shi:1},
we partition $V$ into regions, called $ST$-regions ($ST$ stands for
sign type), as follows.
Let $x \in V$.  For each positive root $\al \in \ro^+$ , 
let $X_{\al}$ equal $+, 0,$ or $ -,$ according to whether
$\langle \al, x \rangle \geq 1$, 
$0 \leq \langle \al, x \rangle < 1$, or $\langle \al, x \rangle <0,$
respectively.
The set of values $\{ X_{\al} \ | \  \al \in  \ro^+ \}$ which arise in
this fashion is called an admissible
sign type.  Now $V$ can be partitioned according to sign type (we drop
the word admissible):  two elements $x, y \in V$ belong to the same
$ST$-region if they determine the same sign type.
It is clear that the interior of an $ST$-region coincides with a 
connected component of the hyperplane arrangement studied by Shi \cite{shi:1}.  

In this note we are interested in the dominant sign types (and 
corresponding $ST$-regions).
These are the sign types with no minus signs; the corresponding regions lie
in $\chamber$.   The dominant sign
types have an important interpretation in terms of the Lie algebra $\g$ of $G$.
Namely, let $\g_{\al} \subset \g$ be the root space for $\al \in \ro$.
Then for each dominant sign type $\sign = \{X_{\al}\}$
the subspace $$\nilrad_{\sign}:=\bigoplus_{X_{\al}=+} \g_{\al}$$ 
is a $B$-stable ideal in the  nilradical $\nilrad$ 
of $\bo$;
and every such ideal arises in this way.
The fact that every ideal arises in this way
can be traced back to Shi's characterization
of admissible sign types \cite{shi:1};  
more recently, it was proved in a different manner by Cellini-Papi 
\cite{cellini-papi:1}.%[proposition].

We now define a partition of $\chamber$ into regions called
$N$-regions indexed by the nilpotent
orbits in $\g$; these are closely related to Lusztig's $\tilde S$-cells in 
the affine Weyl group (see the last section).
Each $N$-region will be a union of certain $ST$-regions. 
Let $\sign$ be a dominant sign type and let $\nilrad_{\sign}$ 
be the corresponding $B$-stable ideal.  Let $\orbit_{\sign}$ denote the
unique nilpotent orbit of $\g$ such that  
$\orbit_{\sign} \cap \nilrad_{\sign}$ is dense in $\nilrad_{\sign}$.
% \cite{lusztig_transf}.
Given a nilpotent orbit $\orbit$ 
we define the $N$-region indexed
by $\orbit$ to be 
$$N_{\orbit}:=\{ x \in \chamber \ | \ \orbit_{\sign_x} = \orbit \}$$
where $\sign_x$ is the sign type determined by $x$.

We now recall the definition of the Dynkin elements of $G$ 
(with respect to $T$ and $B$).  
Given a nilpotent orbit $\orbit$, let $e \in \orbit$ and let $\{e,h,f\}$ be
an $\mathfrak{sl}_2$-triple in $\g$ with $h \in \h = \text{Lie}(T)$.  
Identifying
$V$ with the real subspace of $\h$ spanned by the coroots $\co$, 
we always have $h \in V$ 
($h$, in fact, lies in the integral lattice generated by the coroots;
this follows from the representation theory of $\g$ and 
of $\mathfrak{sl}_2$).  %ref?
By altering our choice of $e$, 
we can ensure that $h \in \chamber$.  Then as noted in the introduction,
$h$ depends only on $\orbit$ and is called the Dynkin element for $\orbit$.
%(or of any element $e$ in $\orbit$). 

Our main result is

\begin{thm} \label{main_theorem}
Let $\orbit$ be a nilpotent orbit in $\g$.  Let $N_{\orbit}$ be the
corresponding $N$-region in $\chamber$.  
Let $x$ be a point of minimal Euclidean length in 
the closure of $N_{\orbit}$.
Then $2 x$ is the Dynkin element for $\orbit$.
In particular, $N_{\orbit}$ is non-empty and $x$ is unique.
\end{thm}

%The paper is organized as follows.  We prove the theorem in the next section.  
%In the fourth section we give the workings of an alternate proof which
%is valid in type $A$ and some other types (but that we expect to work in 
%all types).   Finally in the last section we explain the connection with
%cells in the affine Weyl group and offer some possible generalizations.

\section{Proof}

The theorem follows easily from the next two propositions.
If $h \in \chamber$ is the Dynkin element for the orbit $\orbit$, we
call $$\nilrad_h = \oplus_{i \geq 2} \g_{i}$$ 
the {\bf Dynkin ideal} of $\orbit$,
where $\g_i$ is the $i$-eigenspace for the action of $\mbox{ad}(h)$ 
on $\g$.  This is a $B$-stable ideal in $\nilrad$.

\begin{prop}
Let $h \in \chamber$ be the Dynkin element for the orbit $\orbit$.  Then
$\frac{1}{2} h$ lies in $N_{\orbit}$.
\end{prop}

\begin{proof}
Let $\sign$ be the sign type for $\frac{1}{2} h$.  Then
clearly $\nilrad_{\sign}= \nilrad_h$ from the definitions. 
The result follows from the well-known fact that $\orbit \cap \nilrad_h$
is dense in $\nilrad_h$ \cite{carter}.%specific
\end{proof}

Fix a dominant 
sign type $\sign$ and let $\nilrad_{\sign}$ be the corresponding
ideal.  

\begin{prop} \label{main_prop}
Let $e \in \nilrad_{\sign}$ and let $h \in \chamber$ 
be the Dynkin element of the $G$-orbit through $e$.
Then $\sz{\frac{1}{2} h} \leq \sz{x}$ for any $x$ in the 
closure of the $ST$-region corresponding to $\sign$.  If equality holds,
then $x = \frac{1}{2} h$.
\end{prop}

\begin{proof}
The proof will be a slight 
modification of an argument due to Barbasch and Vogan
(\cite{bv:unip} Lemma 5.7). 
Since all 
$W$-invariant 
positive-definite inner products $( \ , \ )$ on $V$
agree up to a non-zero scalar, we
may as well assume that $(h, \coal)$ is an integer
for all $\coal \in \co$;
this is possible since $h$ lies in the coroot lattice.
Let $\lambda \in \h^*$ 
be the unique element such 
that $\lambda(x) = (h, x)$ for all $x \in V$.
Then $\lambda$ is integral and it is also dominant
since $h \in \chamber$.

Let $(F_{\lambda}, \pi_{\lambda})$ be the irreducible finite-dimensional 
representation of $\g$ of highest weight $\lambda$.
The lowest weight of $F_{\lambda}$ is $w_0( \lambda)$ 
where $w_0$ is the longest element of $W$
(with respect to the simple reflections defined by $\Pi$).
Let $\mathfrak{s} = \{ e', h, f \}$ 
be an $\mathfrak{sl}_2(\complex)$-triple with $e' \in \orbit_e$.
Restricting $F_{\lambda}$ to $\mathfrak{s}$, we see that
the irreducible subrepresentation of $\mathfrak{s}$ of largest dimension 
has dimension $m+1$ where $m = \lambda(h)$ (as $h \in \chamber$).
This us allows to determine the index of nilpotency of $\pi_{\lambda}(e')$ on 
$F_{\lambda}$:  we have 
$$\pi_{\lambda}(e')^m \neq 0, \text{ but } \pi_{\lambda}(e')^{m+1}= 0.$$
Because $e$ and $e'$ are $G$-conjugate, the same holds for $\pi_{\lambda}(e)$.
The key point, to be exploited shortly, is that $m = \lambda(h) = \sz{h}^2$.

On the other hand, let $x$ be in the $ST$-region for $\sign$
(or in fact, its closure).
Viewing $x \in \h$, 
we can filter $F_{\lambda}$ by the eigenspaces of the semisimple operator
$\pi_{\lambda}(x)$ with corresponding eigenvalues greater than a fixed eigenvalue.
The length of this filtration is $\lambda(x) - w_0 \lambda(x)$,
as $x \in \chamber$.

Now if $\al \in \posroots$ 
satisfies $\al(x) \geq 1$, then any non-zero $e_{\al}$ in $\g_{\al}$ acts on 
$F_{\lambda}$ by
raising degrees by at least $1$.
Write $e= \sum e_{\al}$ where $e_{\al} \in \g_{\al}$.
The fact that $e \in \nilrad_{\sign}$ and $x$ belongs
to the closure of the $ST$-region for $\sign$ means that $\al(x) \geq 1$
for all $\al$ appearing in the sum for $e$.
Consequently $e$ acts by raising degrees by at least $1$.
It follows that 
$\pi_{\lambda}( e)^{M+1}=0$ where $M \leq \lambda(x) - w_0 \lambda(x)$. 
The latter quantity is bounded by $2 \sz{h}\sz{x}$ by Cauchy-Schwartz
and the fact that $\sz{w_0(h)} = \sz{h}$.

Hence $m \leq M \leq 2 \sz{h}\sz{x}$.
It follows that $\sz{h}^2 \leq 2\sz{h}\sz{x}$ 
or $\sz{\frac{1}{2} h} \leq \sz{x}$.

If equality holds, then $(h,x) = \sz{h}\sz{x}$ and so $x$ must be a 
positive multiple of $h$.  Hence, $x =  \frac{1}{2} h$.
\end{proof}

The proof of the theorem follows:  every element in the closure of 
$N_{\orbit}$ has size at least the size of one-half the Dynkin element
of $\orbit$ and this minimum is achieved at one-half the Dynkin element.
Uniqueness is a consequence of the 
last statement of the proposition.

We have the following corollary:

\begin{cor}
If $\orbit'$ lies in the closure of the nilpotent orbit $\orbit$
and $\orbit' \neq \orbit$, 
then the size of the Dynkin element of $\orbit'$ 
is strictly smaller than the size of the Dynkin element of $\orbit$.
\end{cor}

\begin{proof}
Let $\nilrad_h = \oplus_{i \geq 2} \g_{i}$ be the Dynkin ideal 
for the Dynkin element $h$ of $\orbit$.
The $G$-saturation of $\nilrad_h$ is by definition the set of elements $\Ad(g)e$ 
where $g \in G$ and $e \in \nilrad_h$;  
this is exactly the closure of $\orbit$. %ref
Hence $\orbit'$ intersects $\nilrad_h$
and the result follows from the previous proposition
applied to the sign type determined by $\frac{1}{2} h$.
\end{proof}

\section{Property D}

In this section, we define a property for $B$-stable ideals that, if it holds in general,
leads to a further understanding of the Dynkin elements.  Furthermore, it should 
also be of interest in the study of $B$-stable ideals and their associated orbits.

Let $\nilrad_{\sign}$ be the ideal corresponding
to the dominant sign-type ${\sign}$.  
Let $\orbit_{\sign}$ be the associated nilpotent orbit;
that is, the one for which $\orbit_{\sign} \cap \nilrad_{\sign}$
is dense in $\nilrad_{\sign}$.

\begin{defn}
We say $\nilrad_{\sign}$ has {\bf property D} if there
exists $e \in \orbit_{\mathcal S} \cap \nilrad_{\mathcal S}$
and $h \in \bo$ such that $e$ and $h$ can be extended to an 
${\mathfrak sl}_2(\complex)$-triple $\{ e, h, f \}$. 
% clarify definitions
%In other words, there exists an ${\mathfrak sl}_2$ triple 
%with a Borel subalgebra contained in $\nilrad_{\mathcal S}$
%such that its nilradical (generically) lies in the $\orbit_{\mathcal S}$.
We say $\nilrad_{\mathcal S}$ has {\bf strong property D} if 
the above holds for any orbit $\orbit$ which intersects $\nilrad_{\mathcal S}$.
\end{defn}

We suspect that every $B$-stable ideal in $\nilrad$ 
possesses strong property D for all $G$.  Note
that it is not true in general that for every element $e \in 
\nilrad_{\mathcal S}$ there exists 
$h \in \bo$ such that  $e$ and $h$ can be extended to an 
${\mathfrak sl}_2(\complex)$-triple.
%counterexample due to McGovern

\begin{lem} \label{lemma1}
Property D holds for every Dynkin ideal.
%and every nilradical of a parabolic subalgebra. % or just distinguished
\end{lem}

\begin{proof}
Let $h \in \h$ be the Dynkin element for the orbit $\orbit$ and 
let $\{e, h , f\}$ be an ${\mathfrak sl}_2(\complex)$-triple.
In particular, $e \in \orbit$.
Then the equation $[h, e] = 2e$ implies $e \in \nilrad_h$,
and thus $e$ and $h$ satisfy the conditions for property D.
\end{proof}

%define basic move
Let $I$ and $I'$ be two $B$-stable ideals in $\nilrad$. 
Suppose that $I = I' \oplus \g_{\beta}$
for some $\beta \in \posroots$. 
For $\al \in \Pi$,
let $P_{\al}$ denote the minimal parabolic subgroup containing $B$
whose Lie algebra contains $\g_{-\al}$.
Assume that $I$ is invariant under the action of $P_{\al}$.
It is not hard to see that $\langle \beta, \coal \rangle \leq 0$.
%If the inequality is strict, we say that $I$ and $I'$ are related by a Texas
%one-step.  

\begin{lem} \label{lemma2}
Assume that $\langle \beta, \coal \rangle < 0$ above.
Then if property D (or the strong version) holds for $I$,
it also holds $I'$.
\end{lem}

\begin{proof}
Let $e \in I$ and $h \in \bo$ satisfy the conditions for
strong property D for the orbit $\orbit_e$.
Since $h$ lies in some Cartan subalgebra of $\bo$, 
we may assume $h \in \h$ since all Cartan subalgebras of $\bo$ are conjugate
under $B$ and $I$ is $B$-stable.

Let $L_{\al}$ be the Levi subgroup of $P_{\al}$ containing $T$.  
Consider the representation of $L_{\al}$ on $I$ and let  
$M$ be the irreducible constituent containing $\g_{\beta}$:
it has lowest weight $\beta$ and highest weight 
$\beta -  \langle \beta, \coal \rangle \al$,
and these are distinct by the hypothesis.

Write $e = \sum e_{\gamma}$, a sum of non-zero root vectors.
If $\gamma = \beta$ does not appear in the sum, we are done 
(since then $e \in I'$).
If $\gamma = \beta -  \langle \beta, \coal \rangle \al$
does not appear in the sum, we proceed by letting 
$n \in N_G(T)$ represent the simple reflection $s_{\al} \in W$.
Then $\Ad(n)e \in I'$ and $\Ad(n)h \in \h$ and we are done.
Finally if both the highest weight and lowest weight for $M$ 
appear in the sum for $e$, we must have $\al(h) = 0$
(as $[h, e] = 2e$)
and thus $h$ lies in the center of the Lie algebra of $L_{\al}$.
From the representation theory of $SL_2(\complex)$
with respect to the representation $M$, there
exists $g \in L_{\al}$ such that $\Ad(g)e \in I'$.
In this case, $\Ad(g)h = h$ and the proof is complete.
\end{proof}

The proof, in particular, implies that $I$ and $I'$ have the same
associated orbit.  
The idea for considering two ideals related in the above fashion
is closely connected to an equivalence relation
on ideals defined in \cite{sommers:equivalences}.

%The three lemmas together with the first 
%theorem in \cite{sommers:equivalences} imply that

\begin{prop}  \label{type-A}
When $G$ is of type $A_n$,
strong property D holds for every $B$-stable ideal in $\nilrad$.
\end{prop}

\begin{proof}
Let $P$ be a parabolic subgroup containing $B$
and $\nilrad_P$ the nilradical of its Lie algebra.  It is not difficult
($G$ is of type $A$) to show that $\nilrad_P$ has strong property D by exhibiting
explicit elements $e \in \orbit$ in a standard form
for each orbit $\orbit$ in the closure of the Richardson orbit for $P$.  
Of course, this uses 
the classification of nilpotent orbits in type $A$ and knowledge of
the closure relations.

Then the first theorem in \cite{sommers:equivalences} 
and the previous lemma imply that strong property D holds for every
$B$-stable ideal in $\nilrad$.   If we wanted only property D, Lemma 
\ref{lemma1} would allow us to avoid the classification of nilpotent orbits.
\end{proof}

We have also checked that strong property D 
holds for all $B$-stable ideals in $\nilrad$ for many $\g$,
including $G_2$, $F_4$, $E_6$, and $E_7$.
We suspect that strong property D always holds.  
One application of property $D$ is the following:

\begin{prop}
Let $\sign$ be a dominant sign type.
If $\nilrad_{\sign}$ possesses property D, 
then there is a convex region containing the closure of the $ST$-region of $\sign$
whose unique minimal point is one-half a Dynkin element for $\orbit_{\sign}$
(here and only here, we do not assume that the Dynkin element 
lies in the closed dominant Weyl chamber).
\end{prop}

\begin{proof}
Let $e \in \nilrad_{\sign}$ and $h \in \bo$ satisfy 
the conditions of property D 
for the orbit $\orbit_e$;  let $f \in \g$ be
such that $\{ e, h, f \}$ is an ${\mathfrak sl}_2(\complex)$-triple. 
We may assume $h \in \h$ as in Lemma \ref{lemma2}
(note that $h$ need not lie in $\mathfrak C$).

Write $e = \sum e_{\al}$ where $e_{\al} \in \g_{\al}$ and
denote by $\supp(e)$ the set of $\al$ which appear in this sum.
Clearly, $\al(\frac{1}{2}h) = 1$ for all $\al \in \supp(e)$.

Define the following convex region in $V$:
$$R_e = \{ v \in V \ | \ \al(v) \geq 1  \text{ for all } \al \in \supp(e) \}.$$
Since $e \in \nilrad_{\sign}$,  $R_e$ contains the closure of the 
$ST$-region for $\sign$.  Then the proof of Proposition  \ref{main_prop}
yields that every point in $R_e$ has length at least the length of $\frac{1}{2}h$.
Moreover, $\frac{1}{2}h$ actually belongs to $R_e$.  Uniqueness also follows from 
Proposition  \ref{main_prop} (or the convexity of $R_e$).
\end{proof}

\section{Relation to cells and representation theory}

Let ${^L \g}$ be the Lie algebra which is of type dual to $\g$.
Let ${^L \h} \subset {^L \g}$ be a Cartan subalgebra.  We identify 
${^L \h}^* \cong \h$.  Let ${^L \bo}\subset {^L \g}$ be the Borel subalgebra
determined by $\bo$.

As noted in the introduction, if $\lambda \in \h$ is one-half a Dynkin element, then 
the weight $\lambda \in {^L \h}^* \cong \h$ is important in representation theory.  
For example, 
let $L(\mu)$ be the irreducible representation of ${^L \g}$ 
of highest weight $\mu \! - \! \rho$, 
where $\rho$ is one-half the sum of the positive roots in ${^L \g}$ determined by ${^L \bo}$.
Then if $\lambda$ is one-half a Dynkin element, the
annihilator $J(\lambda)$ of $L(\lambda)$ is a maximal ideal in ${\mathbf U}({^L \g})$, the universal enveloping
algebra of ${^L \g}$, and this ideal is known to be completely prime in the
classical groups by the work of McGovern (\cite{mcgovern:memoirs} Corollary 6.21). 
%In type $A_n$ this follows from earlier work of Moeglin \cite{moeglin:primitive}.
The idea for studying these particular weights arises in the work of Barbasch and Vogan
on the unitary dual for complex Lie groups;  they are known to yield 
(in many cases, see \cite{bv:unip}) unitary representations.  These representations are 
called special unipotent representations.
In what follows, we locate other weights by varying the type of cell under consideration;  
we hope that these weights might be useful for the unitarity question for
complex Lie groups (and more generally for Arthur's conjectures
for groups over local fields).  
As a first step, one might hope that the purely algebraic work of McGovern extends to 
these weights.

We begin by re-phrasing our main result in type $A_n$.
% by invoking the work of Shi and Lusztig.
Let $\theta$ be the highest root of $\posroots$ and let $s_{0}$ be the affine
reflection in the hyperplane $H_{\theta,1}$ in $V$ where $\theta$ takes the value $1$.
Let $s_1, \dots, s_n$ denote the reflections in  
the hyperplanes $H_{\al, 0}$ where $\al$ takes the value $0$, as $\al$ 
runs through the simple roots in $\Pi$.
Let $\tilde S = \{ s_0, s_1, \dots, s_n \}$.
Then the affine Weyl group of $G$ is the group generated by 
the elements in $\tilde S$; this is a Coxeter group with length function denoted $l(-)$.
For $w \in W_a$, let $R(w) = \{ s \in \tilde S \ | \ l(ws) < l(w) \}$.

Following \cite{lusztig-xi}, we call $\mathcal C \subset W_a$ a canonical left Kazhdan-Lusztig cell
if it is the intersection of a two-sided Kazhdan-Lusztig cell 
and the set of elements $w \in W_a$ such that $R(w) \subset \{s_0\}$.
The canonical left cells are indexed by nilpotent orbits for the group dual to $G$; 
this is 
Lusztig's parametrization of two-sided cells \cite{lusztig:cells-aff-4}.
For example, the cell consisting solely of the identity element in $W_a$
corresponds to the regular nilpotent orbit in ${^L \g}$.

Let $A$ be the fundamental alcove for $W_a$ defined by 
$$A:=\{ x \in V \ | \ \al(x) > 0 \text{ for } \al \in \Pi, \ \theta(x)<1 \}$$
and let $\mathcal C$ be a canonical left Kazhdan-Lusztig 
cell in $W_a$. 
Define 
$\bar{\mathcal C}$ to be the closure of 
$$\{ w(A) \ | \ w \in \mathcal{C} \}$$ 
where the action of $W_a$ on alcoves is the one in \cite{lusztig:cells-aff-1}.
In particular, $\bar{\mathcal C} \subset \chamber$.

By the classification of left Kazhdan-Lusztig cells 
in type $A_n$ 
(due to Shi \cite{shi:book} and Lusztig \cite{lusztig:kac_volume}) 
and by a comment of Lusztig in \cite{lusztig:kac_volume}
(proved by Lawton in \cite{lawton}), 
$\bar{\mathcal C}$ exactly coincides with the closure  of $N_{\orbit}$ where
$\orbit$ is the orbit dual to the one parametrizing $\mathcal C$. 
We remark that the Lawton-Lusztig result is closely related to 
property D.
%and that, in fact, property D can be deduced from it.
% and 
%conversely property $D$ can be used to give another proof of 
%the result from \cite{lawton}.
%(see \cite{sommers:equivalences} for a proof).

Our main theorem now becomes:

\begin{thm}
Let $\orbit$ be a nilpotent orbit in $\mathfrak{sl}_n(\complex)$
and let $\mathcal C$ be the canonical left cell attached to the orbit
which is dual to $\orbit$.  Then the element of minimal length in 
$\bar{\mathcal C}$ is exactly one-half the Dynkin element for $\orbit$.
\end{thm}

Already in $B_2$ and $G_2$
when $\mathcal C$ is a canonical left cell, it is no longer true that 
$\bar{\mathcal C}$ is a union of the closures of $ST$-regions.
Nevertheless, the element $\lambda$ of minimal length in $\bar{\mathcal C}$
is interesting:  the maximal ideal $J(\lambda)$ turns out to be completely prime.  
We can find these weights explicitly in rank $2$ by looking at the pictures
in \cite{lusztig:cells-aff-1}.

\begin{exam}  
In type $B_2$, the elements
are (written as a weighted diagram-- see the introduction -- with the first simple root
being long):  $(0,0)$, $(0, \frac{1}{2})$, $(1, \frac{1}{2})$, $(1,1)$,
coming from the cell attached by Lusztig to 
the regular, the subregular, the minimal, and the zero orbit, respectively.
Only the third element is not one-half a Dynkin element;
the corresponding maximal ideal $J(\lambda)$ is the Joseph ideal,
the unique completely prime primitive ideal in $B_2$
whose associated variety is the minimal orbit \cite{joseph:his_ideal}.
\end{exam}

\begin{exam}  
In type $G_2$, the elements
are (written as a weighted diagram, with the first simple root
being long):  $(0,0)$, $(\frac{1}{2}, 0)$, $(\frac{1}{2}, \frac{1}{2})$, 
$(1, \frac{1}{3})$, $(1,1)$,
coming from the cell attached to 
the regular, the subregular, the eight-dimensional, the minimal, 
and the zero orbit respectively.
The third and fourth elements are not one-half of a Dynkin element.
The fourth element
yields the Joseph ideal \cite{joseph:his_ideal}
and the third element yields a maximal ideal 
which Joseph showed is completely prime and
whose associated variety is the eight-dimensional nilpotent orbit \cite{joseph:2}.% (see also .
\end{exam}

In the rest of the paper, $G$ returns to being of general type.
The above examples suggest 
%that the point of minimal length in 
%a canonical left Kazhdan-Lusztig cell is unique and 

\begin{conj}
Let $\mathcal{C}$ be a canonical left Kazhdan-Lusztig cell in $W_a$.
Let $\lambda \in {^L \h}^* \cong \h$ be an element 
of minimal length in $\bar{\mathcal C}$ (which we conjecture is unique).
Then $J(\lambda)$ is a completely prime, maximal ideal
of ${\mathbf U}({^L \g})$ whose associated variety is
the orbit which indexes $\mathcal{C}$.
\end{conj}

When $G$ is not of type $A_n$, Joseph has shown that there is a unique completely prime, 
primitive ideal whose associated variety is the minimal orbit.
Hence the conjecture predicts that 
when $G$ is not of type $A_n$, the point
of minimal length in the canonical left cell region attached to the minimal orbit
is the infinitesimal character of the Joseph ideal.  
%The values of these characters are listed \cite{joseph:his_ideal}.

We also note that it seems likely that the element of minimal length in $\bar{\mathcal{C}}$
is one-half a Dynkin element 
if and only if the orbit ${^ L} \orbit$ 
indexing $\mathcal{C}$ is special;  in which case,
the minimal element should be one-half the Dynkin element for 
the smallest orbit
in $\g$ which maps under Lusztig-Spaltenstein duality to ${^L}\orbit$.

\medskip

We now investigate another way to get weights which may be of interest in representation
theory.  They may be especially useful for defining unipotent representations 
for complex Lie groups.  It would be interesting to see if there is any 
connection between these weights
and the weights produced in \cite{mcgovern:memoirs} in the classical groups.

In \cite{lusztig:affine_spring}, 
Lusztig defined the notion of $\tilde{S}$-cells in $W_a$.
They are parametrized by pairs $(\orbit, C)$ where $\orbit$
is a nilpotent orbit in $\g$ and $C$ is a conjugacy class in $A(\orbit)$ (
the fundamental group of $\orbit$).  
Let $\mathcal D$ be an $\tilde{S}$-cell 
and let $\bar{\mathcal D} \subset V$ be the closure of 
$$\{ w(A) \ | \ w \in \mathcal{D} \text{ and } w(A) \subset \chamber \}.$$ 
%inverses
We call $\bar{\mathcal D}$ a (canonical left) $\tilde{S}$-region.
It is not too hard to see that the closure of an $N$-region is a union
of  $\tilde{S}$-regions.  More explicitly, if $\orbit$ is a nilpotent orbit in $\g$,
$$\bar{N}_{\orbit} = \bigcup \bar{\mathcal D}_{(\orbit, C)}$$
where the union is over all conjugacy classes $C$ in $A(\orbit)$.
It follows from our main theorem that 
one-half the Dynkin element of $\orbit$ is 
the minimal element of $\bar{\mathcal D}_{(\orbit, C)}$
for some conjugacy class $C \subset A(\orbit)$ (in fact, $C=1$ but we omit
the proof).
This leads us to conjecture 

\begin{conj}
Let $\lambda \in {^L \h}^* \cong \h$ be an element 
of minimal length
in a canonical left $\tilde{S}$-region.
Then $J(\lambda)$ is completely prime.
\end{conj}

In types $B_2$ and $G_2$ the only elements we obtain 
by looking at $\tilde{S}$-regions which are not one-half a Dynkin element
are the elements listed above coming from the Kazhdan-Lusztig cells
(we found this from the $B_2$-picture in \cite{lusztig:affine_spring}
and calculating the canonical left $\tilde{S}$-cells in $G_2$).
This leads us to wonder if the elements from the first conjecture are
a subset of the elements coming from the second conjecture.  This would
follow from a conjecture of the 
second author \cite{sommers:conj_cells} which we now describe. 

Let $\mathcal{N}_{o,c}$ denote the
set of all pairs $(\orbit, C)$ as above.
Let ${^L \mathcal N}_{o}$ denote the nilpotent
orbits for the dual group $^L G$ in ${^L \g}$.
Let $d: \mathcal{N}_{o,c} \to {^L \mathcal N}_{o}$
be the map from \cite{sommers:duality}.

\begin{conj}
Let $\mathcal C$ be the canonical left Kazhdan-Lusztig  
cell indexed by the orbit ${^L \orbit} \subset {^L \g}$.
Then
$$\bar{\mathcal C} = \bigcup \bar{\mathcal{D}}_{(\orbit, C)}$$
%where ${D}_{(\orbit, C)}$ is the canonical left $\tilde{S}$-cell 
%indexed by $(\orbit, C)$ 
where the union is over all pairs $(\orbit , C)$
such that $d_{(\orbit, C)} = {^L \orbit}$.
\end{conj}

\bibliography{sommers_gunn}
\bibliographystyle{pnaplain}
\end{document}